\documentclass[12pt]{amsart}
\usepackage{fullpage,url,amssymb,amsmath,amsfonts,amsthm,mathrsfs}

\usepackage[all]{xy}
\usepackage{color}

\numberwithin{equation}{section}

\renewcommand{\AA}{\mathbb A}
\newcommand{\CC}{\mathbb C}

\newcommand{\FF}{\mathbb F}
\newcommand{\GG}{\mathbb G}

\newcommand{\PP}{\mathbb P}
\newcommand{\QQ}{\mathbb Q}

\newcommand{\ZZ}{\mathbb Z}

\def\bbar#1{\setbox0=\hbox{$#1$}\dimen0=.2\ht0 \kern\dimen0 
\overline{\kern-\dimen0 #1}}
\newcommand{\Qbar}{{\overline{\mathbb Q}}} 
 
\newcommand{\kbar}{\overline{k}} 
\newcommand{\Xbar}{X_{\kbar}}
\newcommand{\FFbar}{\overline{\FF}}  
\newcommand{\KXperp}{K_X^\perp}
\newcommand{\minusK}{\omega_X^{-1}}

\newcommand{\OO}{\mathcal O}
\newcommand{\calE}{\mathcal E}

\newcommand{\sheafO}{\mathscr O}
\newcommand{\Cc}{\mathscr C} 
\newcommand{\F}{\mathscr F}

\newcommand{\Hh}{\mathscr H} 

\newcommand{\Ii}{\mathscr I}

\newcommand{\Ls}{\mathscr L} 
\newcommand{\X}{\mathcal{X}}

\def\NS{\operatorname{NS}} 
\def\Spec{\operatorname{Spec}}

\def\Proj{\operatorname{Proj}}
\def\id{\operatorname{id}}
\def\Div{\operatorname{Div}}
\def\tr{\operatorname{tr}}

\def\Gal{\operatorname{Gal}}
\def\Pic{\operatorname{Pic}}

\def \GL {\operatorname{GL}}

\def\Aut{\operatorname{Aut}} 
\def\End{\operatorname{End}}

\def\Fr{\operatorname{Fr}}

\newcommand{\smallpmod}[1]{\text{ }(\operatorname{mod } #1 )}
\newcommand{\tors}{{\operatorname{tors}}}
\newcommand{\Or}{\operatorname{O}}

\newcommand{\defi}[1]{\textsf{#1}} 

\newtheorem{thm}{Theorem}[section]
\newtheorem{lemma}[thm]{Lemma}

\newtheorem{prop}[thm]{Proposition}

\theoremstyle{definition}
\newtheorem{definition}[thm]{Definition}

\theoremstyle{remark}
\newtheorem{remark}[thm]{Remark}

\newenvironment{romanenum}{\hfill \begin{enumerate} }{\end{enumerate}}

\definecolor{webcolor}{rgb}{0,0,1}
\usepackage[
        colorlinks,
        linkcolor=webcolor,  filecolor=webcolor,  citecolor=webcolor, 
        backref,
        pdfauthor={David Zywina}, pdfauthor={Anthony V\'arilly-Alvarado} 
       pdftitle={Arithmetic $E_8$ lattices with maximal Galois action},
]{hyperref}
\usepackage[alphabetic,backrefs,lite]{amsrefs} 

\renewcommand{\today}{March 19, 2008.}

\begin{document}
\title[]{Arithmetic $E_8$ lattices with maximal Galois action}
\tableofcontents 

\author{Anthony V\'arilly-Alvarado}
\author{David Zywina}
\thanks{The first author was partially supported by a Marie Curie Research Training Network within the 6th European Community Framework Program. The second author was supported by an NSERC postgraduate scholarship.} \address{Department of Mathematics, University of California, 
        Berkeley, CA 94720-3840, USA}
\email{varilly@math.berkeley.edu}
\urladdr{http://math.berkeley.edu/\~{}varilly}
\address{Department of Mathematics, University of California, 
        Berkeley, CA 94720-3840, USA}
\email{zywina@math.berkeley.edu}
\urladdr{http://math.berkeley.edu/\~{}zywina}
\date{\today}

\begin{abstract} 
We construct explicit examples of $E_8$ lattices occurring in arithmetic for which the natural Galois action is equal to the full group of automorphisms of the lattice, i.e., the Weyl group of $E_8$.  In particular, we give explicit elliptic curves over $\mathbb{Q}(t)$ whose Mordell-Weil lattices are isomorphic to $E_8$ and have maximal Galois action.

Our main objects of study are del Pezzo surfaces of degree $1$ over number fields. The geometric Picard group, considered as a lattice via the negative of the intersection pairing, contains a sublattice isomorphic to $E_8$. We construct examples of such surfaces for which the action of Galois on the geometric Picard group is maximal. 
\end{abstract}

\maketitle

\section{Introduction} \label{S:intro}
In this paper we construct explicit examples of $E_8$ lattices coming
from arithmetic for which the Galois action is as large as possible.
Recall that a \defi{lattice} is a free abelian group equipped with a $
\ZZ$-valued non-degenerate symmetric bilinear form.  The $E_8$
lattice is the unique positive definite, even, unimodular lattice of
rank $8$.   As the notation suggests, it is also the root lattice of
the $E_8$ root system (which is the largest exceptional root
system).  The \defi{Weyl group} of $E_8$, denoted $W(E_8)$, is the
group of automorphisms of the $E_8$ lattice; it is a finite group of
order $696,\!729,\!600$.

Suppose that $E$ is a
non-isotrivial elliptic curve over $\QQ(t)$.  The group
$E(\Qbar(t))/E(\Qbar(t))_\tors$ is free abelian of finite rank, and
has a natural pairing called the \defi{canonical height pairing} (see~\cite[Theorem III.4.3]{Silverman2}).
The group $E(\Qbar(t))/E(\Qbar(t))_\tors$ equipped with this pairing
is a lattice, called the \defi{Mordell-Weil lattice} of $E$. The natural action
of $\Gal(\Qbar/\QQ)$ on $E(\Qbar(t))/E(\Qbar(t))_\tors$ preserves the lattice structure.

Let us now give an example of an $E_8$ lattice occurring in
arithmetic.  

\begin{thm} \label{T:E8}
Let $a(t), b(t), c(t) \in \ZZ[t]$ be polynomials of degrees at
most $2,4,$ and $6$, respectively, which satisfy the following
congruences:
\begin{align*}
a(t) &\equiv 70t \pmod{105}, \\
b(t) &\equiv 84t^4 + 7t^3 + 98t^2  + 84t + 98 \pmod{105}, \\
c(t) &\equiv 65t^6 + 42t^5 + 21t^4 + 77t^3 + 63t^2 + 56t + 30 \pmod
{105}.
\end{align*}
Let $E$ be the elliptic curve over $\QQ(t)$ given by the Weierstrass
model
\begin{equation} \label{E:W eqn for E8}
y^2 = x^3 + a(t)x^2 + b(t)x + c(t).
\end{equation}
Then the group $E(\Qbar(t))$ is free of rank $8$ and as a
Mordell-Weil lattice it is isomorphic to the $E_8$ lattice.   The group
$\Gal(\Qbar/\QQ)$ acts on $E(\Qbar(t))$ as the full group $W(E_8)$.
\label{T:curvewithsymmetry}
\end{thm}

\begin{remark}
\begin{romanenum}
\item 
Mordell-Weil lattices of type $E_8$ have been extensively studied by Shioda (see \cite{Shioda, Shioda3}).  In fact, Shioda \cite{Shioda}*{Theorem~7.2} has shown that \emph{every} Galois extension of $\QQ$ with Galois group isomorphic to $W(E_8)$ arises from the Mordell-Weil lattice of an elliptic curve $E/\QQ(t)$.  Theorem~\ref{T:E8} thus gives explicit examples of this theory.
\item
If the Mordell-Weil lattice of $E/\QQ(t)$ is isomorphic to the $E_8$ lattice, then the $240$ roots of the lattice are $\Qbar(t)$-rational points of the form
\begin{equation}
\label{E:points}
x = a_2t^2 + a_1t + a_0 \qquad y = b_3t^3 + b_2t^2 + b_1t + b_0,
\end{equation}
where $a_i, b_j \in \Qbar$ for all $i$ and $j$. Conversely, any $\Qbar(t)$-point of this form is a root of the lattice; see~\cite[\S10]{Shioda2}. 
The field extension of $\QQ$ obtained by adjoining all the coefficients $a_i$ and $b_j$ is a Galois extension of $\QQ$ with Galois group isomorphic to a subgroup of $W(E_8)$. 

Theorem~\ref{T:curvewithsymmetry} can thus be used, in principle, to write down explicit Galois extensions of $\QQ$ with Galois group isomorphic to $W(E_8)$: first substitute the expressions of~(\ref{E:points}) into~(\ref{E:W eqn for E8}). We then obtain a polynomial in the variable $t$, which is identically zero if $(x,y) \in E(\Qbar(t))$. This will give a series of relations among the $a_i$ and $b_j$. We can then use elimination theory to distill these relations to a single polynomial in, say, the variable $a_1$. The Galois group of the splitting field of this polynomial will be isomorphic to $W(E_8)$.

Another method for obtaining a polynomial whose splitting fields has Galois group isomorphic to $W(E_8)$, can be found in \cite{JKZ}.  This approach uses an algebraic group of $E_8$ type.
\end{romanenum}
\end{remark}

\subsection{Del Pezzo Surfaces}

A \defi{del Pezzo surface $X$ of degree 1} over a number field $k$ is
a surface that when base extended to an algebraic closure $\kbar$ of
$k$ is isomorphic to $\PP^2_{\kbar}$ blown up at $8$ points in general
position.

An $E_8$ lattice naturally arises from a del Pezzo surface $X$ of
degree $1$ as follows. Let $X_{\kbar} = X \times_k \kbar$. The intersection pairing on $\Pic(\Xbar) \cong
\ZZ^9$ gives this group a lattice structure. Let $\KXperp$ denote the
orthogonal complement in $\Pic(\Xbar)$ of the canonical class $K_X$
with respect to the intersection pairing. We give $\KXperp$ a lattice
structure by endowing it with the form that is the {negative} of
the intersection pairing.  As a lattice, $\KXperp$ is isomorphic to
the $E_8$ lattice.  We denote by $\Or(\KXperp)$ the group of lattice
automorphisms of $\KXperp$. It is thus isomorphic to $W(E_8)$.

There is a natural action of $\Gal(\kbar/k)$ on $\Pic(\Xbar)$ that
respects the intersection pairing and fixes $K_X$.  We thus obtain a
Galois representation
\begin{equation*}
\phi_X\colon\Gal(\kbar/k) \to \Or(\KXperp).
\end{equation*}
We will construct explicit examples of $X$ such that $\phi_X$ is
surjective, that is, such that the action of Galois on $\KXperp$ is
maximal. Ekedahl~\cite{Ekedahl} and Ern\'e~\cite{Erne} have made the analogous constructions for del Pezzo surfaces of degree $3$ and $2$, respectively (for the general definition of del Pezzo
surfaces, see \S\ref{S:background}). The reader will see their influence here. In~\cite{Shioda}, Shioda also constructs explicit examples for del Pezzo surfaces of degrees $2$ and $3$, using the theory of Mordell-Weil lattices.

We can now state our main result.
\begin{thm}
\label{T:main theorem}
Let $f$ be a sextic polynomial in the weighted graded ring
$\ZZ[x,y,z,w]$, where the variables $x,y,z,w$ have weights $1,1,2,3$
respectively, such that
\begin{align} \label{E:main equation}
f\equiv 40x^6 &+ 63x^5y + 84x^4y^2 + 28x^3y^3 + 42x^2y^4 + 49xy^5 +
75y^6 + 84x^4z  \\
\notag &+ 7x^3yz
   + 98x^2y^2z  + 84xy^3z + 98y^4z + 35xyz^2 + z^3 + w^2  \pmod{105}.
\end{align}
Then $X := \Proj\left(\QQ[x,y,z,w]/(f)\right)$ is a del Pezzo surface
of degree $1$ over $\QQ$ and the image of the homomorphism
\[
\phi_X\colon \Gal(\Qbar/\QQ) \to \Or(\KXperp)
\]
is surjective.  Equivalently, if $L_X$ is the fixed field of
$\ker(\phi_X)$ in $\Qbar$, then
\[
\Gal(L_X/\QQ) \cong W(E_8).
\]
Letting the polynomials $f$ that satisfy (\ref{E:main equation}) vary,
the corresponding fields $L_X$ give infinitely many linearly disjoint
extensions of $\QQ$ with Galois group isomorphic to $W(E_8)$.
\end{thm}

\begin{remark}
\begin{romanenum}
\item
Suppose that $k$ is a finitely generated extension of $\QQ$.   By Theorem~\ref{T:main theorem}, there is a del Pezzo surface $X/\QQ$ of degree $1$ such that $\Gal(L_X/\QQ)\cong W(E_8)$, and such that $k$ and $L_X$ contain no isomorphic subfields except $\QQ$ (the proof of Theorem~\ref{T:main theorem} will give a constructive way to find such an $X$).  Therefore $X_k$ is a del Pezzo surface of degree $1$ over $k$, and the homomorphism $\phi_{X_k}\colon \Gal(\kbar/k) \to \Or(K_{X_k}^\perp)$ is surjective.
\item
The finite subgroups of $\GL_8(\QQ)$ with maximal cardinality are isomorphic to $W(E_8)$ (this depends on unpublished results in group theory, see \cite{BDEPS}).  Thus our examples have maximal Galois action amongst all rank $8$ lattices.
\end{romanenum}
\end{remark}

\subsection{Genus $4$ curves}
We shall now describe certain genus $4$ curves that arise from del Pezzo surfaces of degree $1$.  These curves have been studied by Zarhin, and our examples complement his.  For details, see \cite{Zarhin}.

Let $X$ be a del Pezzo surface of degree $1$ over a field $k$ of characteristic $0$.  The surface $X$ has a distinguished involution called the \defi{Bertini involution}; it is the unique automorphism of $X$ which induces an action of $-I$ on $\KXperp\subseteq \Pic(\Xbar)$.  The fixed locus of the Bertini involution consists of a curve $C$ and a rational point.
The curve $C$ is smooth, irreducible, non-hyperelliptic, and has genus $4$. 

 Let $J(C)$ be the Jacobian of $C$, and let $J(C)[2]$ be the $2$-torsion subgroup of $J(C)(\kbar)$.  The group $J(C)[2]$ is an $8$-dimensional vector space over $\FF_2$ and is equipped with the \emph{Weil pairing}
\[
\langle\ , \, \rangle \colon J(C)[2]\times J(C)[2] \to \{\pm 1\} \cong \FF_2.
\]
The pairing $\langle\ , \, \rangle$ is an alternating nondegenerate bilinear form.  The Galois group $\Gal(\kbar/k)$ acts on $J(C)[2]$ and preserves the Weil pairing.  The following lemma describes the structure of $J(C)[2]$.

\begin{lemma}[\cite{Zarhin}*{Theorem~2.10 and Remark~2.12}]
\label{L:Zarhin}
There is an isomorphism of $\Gal(\kbar/k)$-modules
\[
J(C)[2] \cong \KXperp/2\KXperp,
\]
which preserves the corresponding $\FF_2$-valued pairings.
\end{lemma}

Recall that $\KXperp \cong E_8$,  so the Galois action on $J(C)[2]$ factors through the group $W(E_8)/\{\pm I\}$.  We may thus use our examples from Theorem~\ref{T:main theorem} to give examples of such curves with maximal Galois action.

\begin{prop}  \label{T:genus 4 curves}
Let $f$ be a sextic polynomial in the weighted graded ring $\ZZ[x,y,z]$, where the variables $x,y,z$ have weights $1,1,2$ respectively, such that 
\begin{align*}
f\equiv 40x^6 + 63x^5y &+ 84x^4y^2 + 28x^3y^3 + 42x^2y^4 + 49xy^5 + 75y^6 + 84x^4z  \\
\notag &+ 7x^3yz  
    + 98x^2y^2z  + 84xy^3z + 98y^4z + 35xyz^2 + z^3  \pmod{105}.
\end{align*}
Define the curve $C := \Proj(\QQ[x,y,z]/(f))$, and let $J(C)$ be its Jacobian.  The curve $C$  is smooth, geometrically irreducible, non-hyperelliptic, and has genus $4$.   There is an isometry
\[
J(C)[2] \cong E_8/2E_8
\]
under which the absolute Galois group $\Gal(\Qbar/\QQ)$ acts as the full group $W(E_8)/\{\pm I\}.$   

The Jacobian $J(C)$ has no non-trivial endomorphisms over $\Qbar$, i.e., $\End(J(C)_\Qbar)= \ZZ$.  In particular, $J(C)$ is an absolutely simple abelian four-fold.
\end{prop}
\begin{proof}
For the del Pezzo surfaces in Theorem~\ref{T:main theorem}, the automorphism $[x:y:z:w] \mapsto [x:y:z:-w]$ is the Bertini involution.  The curves in the proposition are in the fixed locus of the Bertini involution of the appropriate del Pezzo surface.  That the Galois action is maximal, is then a consequence of Theorem~\ref{T:main theorem} and Lemma~\ref{L:Zarhin}.  
The last statements of the Proposition follows from \cite{Zarhin}*{Theorem~4.3}.
\end{proof}

\subsection{Overview}

Let us briefly outline the contents of this paper.   

In \S\ref{S:background}, we summarize some of the basic theory of del Pezzo surfaces, focussing on the aspects relevant to our application.

In \S\ref{S:grouptheory}, we prove a useful criterion to determine whether a subgroup $H\subseteq W(E_8)$ is actually the full group $W(E_8)$.  Applied to a del Pezzo surface $X/\QQ$ of degree $1$, it is gives a criterion for the representation $\phi_X$ to be surjective (Proposition~\ref{P:full Galois criterion}).  In particular, to show that $\phi_X$ is surjective, it suffices to give a model $\mathcal{X}/\ZZ$ of $X$ with certain kinds of reduction at three special fibers.   For example, one of the conditions is the existence of a prime $p$ for which $\X_{\FF_p}$ is a del Pezzo surface of degree $1$ and $\phi_{\X_{\FF_p}} \colon \Gal(\FFbar_p/\FF_p) \to \Or(K^\perp_{\X_{\FF_p}})$ has image of order $7$.
 
In order to construct the surfaces of Theorem~\ref{T:main theorem}, we first construct three del Pezzo surfaces of degree $1$ over the finite fields $\FF_3$, $\FF_5$ and $\FF_7$, with the properties required by Proposition~\ref{P:full Galois criterion} (this explains the congruences modulo $105=3\cdot 5 \cdot 7$).   To prove Theorem~\ref{T:main theorem}, we then exhibit a scheme $\X/\ZZ$ whose generic fiber is $X/\QQ$, and whose special fibers at $3,5,$ and $7$ are isomorphic to the del Pezzo surfaces already calculated (see \S\ref{S:main proof}).  
 
Our surfaces over $\FF_3$ and $\FF_5$ will be given explicitly as the blow-up of $8$ points in the projective plane (see \S\ref{Ex:F3 example} and \S\ref{Ex:F5 example}).   In \S\ref{S:dictionary}, we describe how given the blow-up data, one can write down an explicit (weighted) sextic polynomial defining the surfaces; i.e., the anticanonical model.

For our example over $\FF_7$, we simply write down a candidate surface and then verify that it satisfies the required properties by applying the Lefschetz trace formula (see \S\ref{SS:F7 example}).

Finally in \S\ref{S:Elliptic curve proofs}, we show how given a del Pezzo surface $X/\QQ$ of degree $1$, one can associate an elliptic curve $E/\QQ(t)$.  As a consequence of the work of Shioda, there is a lattice isomorphism $\KXperp \cong E(\Qbar(t))/E(\Qbar(t))_\tors$ which respects the corresponding Galois actions.   Working this out explicitly, we will find that Theorem~\ref{T:E8} is a direct consequence of Theorem~\ref{T:main theorem}.

\subsection*{Acknowledgements} 
We thank Bjorn Poonen for many helpful comments.  Our computations were performed using {\tt Magma}~\cite{magma}.

\subsection*{Notation}
We now fix some notation and conventions which will hold throughout the paper.  For a field $k$, fix an algebraic closure $\kbar$.  For $S$-schemes $X$ and $Y$, define $X_Y:=X\times_S Y$; if $Y=\Spec B$, then we will write $X_B$ for $X_{\Spec{B}}$.    If $\F$ is a sheaf  of $\sheafO_X$-modules on a $k$-scheme $X$, then the dimension of the $k$-vector space $H^0(X,\F)$ of global sections will be denoted  by $h^0(X,\F)$.  

By a \defi{surface}, we mean a smooth projective geometrically integral scheme of dimension $2$ of finite type over a field $k$.  Given a surface $X$ over $k$, we have an intersection pairing $(\ ,\, )\colon \Pic(X_{\kbar})\times\Pic(X_{\kbar}) \to \ZZ.$ We write $\omega_X$ for the canonical sheaf of $X$, and $K_X$ for its class in the Picard group. We will identify $\Pic(X)$ with the Weil divisor class group; in particular, we will use additive notation for the group law on $\Pic(X)$.

Suppose that $k$ is a number field and that $v$ is a finite place of $k$.  Denote the completion of $k$ at $v$ by $k_v$, the valuation ring of $k_v$ by $\OO_v$, and the corresponding residue field by $\FF_v$.  Let $F_v\in \Gal(\FFbar_v/\FF_v)$ be the Frobenius automorphism $x\mapsto x^{|\FF_v|}$.  We will denote the ring of integers of $k$ by $\OO_k$. If $S$ is a set of places of $k$, we write $\OO_{k,S}$ for the ring of $S$-integers of $k$.

\section{Background on del Pezzo surfaces}
\label{S:background}

\subsection{Del Pezzo surfaces}
We now review some basic theory concerning del Pezzo surfaces. The standard references on the subject are~\cite{Manin}, ~\cite{Demazure1980} and~\cite[III.3]{Kollar1996}.  

\begin{definition}
A \defi{del Pezzo surface} over a field $k$ is a surface $X$ over $k$ with ample anticanonical sheaf $\minusK$.  
The \defi{degree} of $X$ is the intersection number $(K_X,K_X)$.
\end{definition}
For a del Pezzo surface $X$ over $k$ of degree $d$, we have $1\leq d \leq 9$.  The surface $\Xbar$ is isomorphic to either $\PP^1_{\kbar}\times\PP^1_{\kbar}$ (which has degree $8$) or to the blow-up of $\PP^2_{\kbar}$ at $r := 9 - d$ closed points. Moreover, in the second case, the $r$ points are in \defi{general position}; that is, no $3$ of them lie on a line, no $6$ of them lie on a conic, and no $8$ of them lie on a cubic with a singularity at one of the points. 

For $r\leq 8$, the blow-up of $r$ distinct closed points of  $\PP^2_{\kbar}$ in general position is a del Pezzo surface of degree $9-r$ over $\kbar$~\cite[Theorem 1, p.27]{Demazure1980}. 
\subsection{Structure of the Picard group}
\label{SS:PicardGroup}

Let $X$ be a del Pezzo surface over $k$ of degree $d \leq 6$.  The Picard group $\Pic(\Xbar)$ is a free abelian group of rank $10 - d = r+1$, and has a basis $e_1,\dots,e_r,\ell$ such that
\[
(e_i,e_j) = -\delta_{ij},\quad 
(e_i,\ell) = 0,\quad
(\ell,\ell) = 1,\quad\text{and \ }
{-K}_X = 3\ell - \sum_{i=1}^re_i.
\]
If $X_{\kbar}$ is the blow-up of $\PP^2_{\kbar}$ along a set of closed points $\{P_1,\dots,P_r\}$ which are in general position, then we may take $e_i$ to be the exceptional divisor corresponding to $P_i$ and $\ell$ to be the strict transform of a line in $\PP^2_{\kbar}$ not passing through any of the $P_i$. The divisors $e_i$ are isomorphic to $\PP^1_{\kbar}$.

\begin{definition} 
Let $\KXperp$ be the orthogonal complement of $K_X$ in $\Pic(\Xbar)$ with respect to the intersection pairing.  We give $\KXperp$ the structure of a lattice by endowing it with the form that is the \emph{negative} of the intersection pairing.
\end{definition}

The lattice $\KXperp$  is isomorphic to the root lattice of type $A_1\times A_2$, $A_4$, $D_5$, $E_6$, $E_7$ or $E_8$ (where $r$ is the sum of the subscripts) \cite[Theorem 23.9]{Manin}.  The group $\Or(\KXperp)$ of lattice automorphisms of $\KXperp$ is isomorphic to the Weyl group of $\KXperp$ (see  \cite[Theorem 23.9 and \S26.5]{Manin}).

\subsection{Galois action on the Picard group}
\label{SS:GaloisAction}
Let $X$ be a del Pezzo surface over $k$. 
For each $\sigma \in \Gal(\kbar/k)$, let $\tilde{\sigma}\colon\Spec\kbar \to \Spec\kbar$ be the corresponding morphism. Then $\id_X\times \tilde{\sigma}\in \Aut(\Xbar)$ induces an automorphism $(\id_X\times \tilde{\sigma})^*$ of $\Pic(\Xbar)$. This action of $\Gal(\kbar/k)$ on $\Pic(\Xbar)$ fixes the canonical class $K_X$ and preserves the intersection pairing.  Therefore it factors through the action on $\KXperp$, inducing a group homomorphism
\begin{align*}
\phi_X\colon \Gal(\kbar/k) &\to \Or(\KXperp) \\
\sigma &\mapsto (\id_X\times \tilde{\sigma})^*|_{\KXperp}.
\end{align*}

\subsection{Weighted projective spaces} \label{SS:wps}
We quickly recall some basic definitions for weighted projective spaces; a good general  reference is \cite{Dolgachev}.
Fix a field $k$ and positive integers $q_0,\dots,q_n$.  Let $k[x_0,\dots,x_n]$ be the polynomial ring in the variables $x_0,\dots,x_n$ graded with weights $q_0,\dots,q_n,$ respectively. We define the weighted projective space $\PP_k(q_0,\dots,q_n)$ to be the $k$-scheme $\Proj(k[x_0,\dots,x_n])$. 

The $\ZZ$-grading on the variables $x_0,\dots,x_n$ gives rise to an action of $\GG_m$ on $\AA^{n+1}_k = \Spec(k[x_0,\dots,x_n])$ via the $k$-algebra homomorphism
\begin{align*}
k[x_0,\dots,x_n] &\to k[x_0,\dots,x_n] \otimes k[t,t^{-1}] \\
x_i &\mapsto x_i \otimes t^{q_i}.
\end{align*}
The open subscheme $U = \AA^{n+1}_k \setminus \{0\}$ is stable under this action.  The universal geometric quotient $U/\GG_m$ exists and coincides with $\PP_k(q_0,\dots,q_n)$.  The set of $\kbar$-valued points of $\PP_k(q_0,\dots,q_n)$ can be described as
\[
\PP_k(q_0,\dots,q_n)(\kbar) = \big( \kbar^{n+1}-\{(0,\dots,0)\} \big) / \sim,
\]
where $(a_0,\dots,a_n) \sim (a'_0,\dots,a'_n)$ if there exists a $\lambda \in \kbar^\times$ such that $a_i = \lambda^{q_i}a'_i$ for all $i$. We denote the equivalence class of $(a_0,\dots,a_n)$ by $[a_0:\ldots:a_n]$.  
A homogeneous ideal $I$ of $k[x_0,\dots,x_n]$ (with respect to the above grading) determines a closed subscheme $V(I) := \Proj(k[x_0,\dots,x_n]/I)$ of $\PP_k(q_0,\dots,q_n)$. 
The set of $\kbar$-valued points of $V(I)$ is
\[
V(I)(\kbar) = \{[a_0:\ldots:a_n] \in \PP_k(q_0,\dots,q_n)(\kbar) \,:\, f(a_0,\dots,a_n) = 0\text{ for all homogeneous } f \in I \}.
\]

\subsection{The anticanonical model}
Besides the blow-up description, there is another useful model of a del Pezzo surface.  For any $k$-scheme $X$ and line bundle $\Ls$ on $X$, we may construct the graded $k$-algebra
\[
R(X,\Ls) := \bigoplus_{m\geq 0} H^0(X,\Ls^{\otimes m}).
\]
When $\Ls = \minusK$, we call $R(X,\minusK)$ the \defi{anticanonical ring} of $X$.  If $X$ is a del Pezzo surface over $k$, then $X \cong \Proj R(X,\minusK)$ \cite[Theorem III.3.5]{Kollar1996}.  The $k$-scheme $\Proj R(X,\minusK)$ is known as the \defi{anticanonical model} of $X$. 

\begin{prop} \label{P:sextic}
Let $X$ be a del Pezzo surface of degree $1$ over a field $k$, and let  $x,y,z,w$ be variables with weights $1,1,2,3$ respectively.  Then there is an isomorphism of graded $k$-algebras
\[
R(X,\minusK)\cong k[x,y,z,w]/(f),
\] 
where $f$ is a sextic in $k[x,y,z,w]$.  The surface $X$ is thus isomorphic to the smooth sextic hypersurface $V(f)$ in $\PP_k(1,1,2,3)$.  

Conversely, if $f \in k[x,y,z,w]$ is a sextic polynomial such that $V(f)$ is smooth, then $V(f)$ is a del Pezzo surface of degree $1$.
\end{prop}
\begin{proof}
See~\cite[Theorem III.3.5]{Kollar1996}.
\end{proof}
In \S\ref{S:dictionary} we will show how, given a blow-up model of a del Pezzo surface $X$ of degree $1$, one can find a sextic polynomial $f$ as in Proposition~\ref{P:sextic}.

\begin{lemma} \label{L:anticanonical map}
Let $X$ be a del Pezzo surface of degree $1$ defined over a field $k$.  The linear system $|{-K}_X|$ has a single base point $O$, and this point is defined over $k$.  We call $O$ the \defi{anticanonical point} of $X$.  
If $X$ is the locus of a sextic $f(x,y,z,w)$ in $\PP_k(1,1,2,3)$, then the the linear system $|{-K}_X|$ gives rise to the rational map
\[
X  \dashrightarrow \PP^1_k, \,[x:y:z:w] \mapsto [x:y].
\]
\end{lemma}
\begin{proof}
For the first statement, see~\cite[p. 40]{Demazure1980}. If $X$ is a sextic $f(x,y,z,w) = 0$ in $\PP_k(1,1,2,3)$ then the functions $x$ and $y$ form a basis for $H^0(X,\omega_X^{-1})$ (see~\cite[Proof of Theorem~3.36(6)]{Kollar}).
\end{proof}
\subsection{Reductions of Del Pezzo surfaces}

\begin{lemma}  \label{L:specialization}
Let $v$ be a finite place of a number field $k$.  Suppose that $\X$ is a smooth $\OO_v$-scheme such that $X:=\X_{k_v}$ and $\X_{\FF_v}$ are both del Pezzo surfaces of the same degree $d\leq 6$.
Given $\sigma \in \Gal(\kbar_v/k_v)$, restrict $\sigma$ to the maximal unramified extension of $k_v$ in $\kbar_v$ and let $\overline{\sigma}\in\Gal(\FFbar_v/\FF_v)$ be the corresponding automorphism of residue fields.  
There exists an isomorphism of lattices
\[
\beta \colon K_{X}^\perp \stackrel{\sim}{\to} 
K_{\X_{\FF_v}}^\perp
\]
such that for any $\sigma \in \Gal(\kbar_v/k_v)$ we have
\[
\phi_{X}(\sigma) = \beta^{-1} \phi_{\X_{\FF_v}}(\overline{\sigma}) \beta.
\]
\end{lemma}
\begin{proof}
Fix an integral divisor $V$ of $X$.  Let $\overline{V}$ be the Zariski closure of $V$ in $\X$.   Since $\overline{V}$ is an $\OO_v$-scheme, we may consider its reduction $\overline{V}_{\FF_v}$.
Extending $V\mapsto \overline{V}_{\FF_v}$ by additivity defines a group homomorphism $\alpha\colon \Pic(X) \to \Pic(\X_{\FF_v})$ called the \defi{specialization map} (see~\cite[\S 20.3]{Fulton}).  By \cite{Fulton}*{Corollary 20.3}, $\alpha$ preserves intersection pairings.  From our description of the intersection pairing on $\Pic(X_{\kbar_v})$ in \S\ref{SS:PicardGroup}, we find that the intersection pairing is nondegenerate and hence $\alpha$ is injective.   

We claim that $\alpha(K_X) = K_{\X_{\FF_v}}$.  There is a commutative diagram \cite[\S 20.3.1]{Fulton},
\[
\xymatrix{
 & \Pic(\X) \ar[dl]_{j^*}\ar[dr]^{i^*} & \\
 \Pic(X) \ar[rr]^{\alpha} & & \Pic(\X_{\FF_v}),
}
\]
where $j\colon X \to \X$ and $i\colon \X_{\FF_v} \to \X$ are the morphisms coming from the respective fiber products.  Thus it is enough to show that $j^*(K_{\X}) = K_X$ and $i^*(K_{\X}) = K_{\X_{\FF_v}}$.  Since pull-backs commute with tensor operations~\cite[Ex. II.5.16(e)]{Hartshorne}, it suffices to prove that 
\[
\Omega^1_{\X/\OO_v}\times_{\OO_v}k_v \cong \Omega^1_{X/k_v}\quad\text{and}\quad
\Omega^1_{\X/\OO_v}\times_{\OO_v}\FF_v \cong \Omega^1_{\X_{\FF_v}/\FF_v};
\]
these isomorphisms follow from the compatibility of relative differentials with base extension~\cite[II.8.10]{Hartshorne}. Therefore $\alpha(K_X)=K_{\X_{\FF_v}}$.   

Since $\alpha$ preserves intersection pairings and $\alpha(K_X)=K_{\X_{\FF_v}}$, we have an injection of lattices
\[
\beta := \alpha|_{\KXperp} \colon \KXperp {\to} K_{\X_{\FF_v}}^\perp.
\]
Both lattices are unimodular,  so this map must also be surjective.

Take any $\sigma \in \Gal(\kbar_v/k_v)$.  For an integral divisor $V$ of $X$, one finds that
$\overline{\sigma(V)}_{\FF_v} = \overline{\sigma}(\overline{V}_{\FF_v});$
equivalently, $\alpha$ commutes with the respective Galois actions.  It is then an immediate consequence that $\beta \phi_{X}(\sigma) =\phi_{\X_{\FF_v}}(\overline{\sigma}) \beta.$
\end{proof}

\begin{lemma} \label{L:GaloisLift}
Let $X$ be a del Pezzo surface of degree $d\leq 6$ over a number field $k$.  Let $S$ be a finite set of places of $k$. Let $\X$ be a smooth $\OO_{k,S}$-scheme for which $X=\X_{k}$, and let $v \notin S$ be a finite place of $k$ such that $\X_{\FF_v}$ is also a del Pezzo surface of degree $d$. 
Then there is a lattice isomorphism $\theta \colon \KXperp \overset{\sim}{\to} K_{\X_{\FF_v}}^\perp$ and an automorphism $\sigma \in \Gal(\kbar/k)$ such that 
\[
\phi_X(\sigma) = \theta^{-1} \phi_{\X_{\FF_v}}(F_v) \theta.
\]
\end{lemma}
\begin{proof}
In the notation of Lemma~\ref{L:specialization}, choose $\sigma \in \Gal(\kbar_v/k_v)$ with $\overline{\sigma}=F_v$.  Applying Lemma~\ref{L:specialization} with $\X_{\OO_v}$, we know there is a lattice isomorphism $\beta \colon K_{X_{k_v}}^\perp \overset{\sim}{\to }
K_{\X_{\FF_v}}^\perp$ for which
\[
\phi_{X_{k_v}}(\sigma) = \beta^{-1} \phi_{\X_{\FF_v}}(F_v) \beta.
\]
Now fix an embedding $\iota\colon\kbar\hookrightarrow \kbar_v$.  This gives an inclusion $\Gal(\kbar_v/k_v) \hookrightarrow \Gal(\kbar/k)$, and hence we may also view $\sigma$ as an automorphism of $\kbar$.  The embedding $\iota$ induces an isomorphism $\gamma\colon \Pic(X_{\kbar}) \overset{\sim}{\to} \Pic(X_{\kbar_v})$ which preserves the intersection pairing and the canonical class.  Therefore $\gamma|_{\KXperp}$ is a lattice isomorphism from $\KXperp$ to $K_{X_{k_v}}^\perp$ such that
\[
\phi_X(\sigma) = (\gamma|_{\KXperp})^{-1} \phi_{X_{k_v}}(\sigma) \gamma|_{\KXperp}.
\]
The lemma follows by setting $\theta:= \beta \gamma|_{\KXperp}$.
\end{proof}

\section{Group theory}
\label{S:grouptheory}

We wish to find a del Pezzo surface $X$ of degree $1$ with surjective homomorphism 
\[
\phi_X\colon\Gal(\kbar/k)\to \Or(\KXperp).
\]  
To accomplish this, we will need a convenient criterion to determine whether a subgroup of $\Or(\KXperp)\cong W(E_8)$ is the full group.\\

The Weyl group $W(E_8)$ may be viewed as a subgroup of $\GL(E_8) \cong \GL_8(\ZZ)$, so we can talk about the trace, determinant, and characteristic polynomials of elements (or conjugacy classes) of $W(E_8)$.  The \emph{order} of a conjugacy class is the order of any element in the conjugacy class as a group element; this should not be confused with the \emph{cardinality} of a conjugacy class which is the number of elements it contains.  The goal of this section is to prove the following group theoretic proposition.
\begin{prop} \label{P:Weyl group criterion}
Let $H$ be a subgroup of $W(E_8)$.  Suppose that the following conditions hold:
\begin{romanenum}
\item \label{I:order 7} There exists an element in $H$ of order $7$,
\item \label{I:trace 5} There exists an element in $H$ of order $3$ and trace $5$,
\item \label{I:trace -4} There exists an element in $H$ of order $3$ and trace $-4$,
\item \label{I:negative determinant} There exists an element in $H$ with determinant $-1$.
\end{romanenum}
Then $H=W(E_8)$.
\end{prop}

In terms of our application to del Pezzo surfaces, we have the following criterion. 
\begin{prop} \label{P:full Galois criterion}
Let $X$ be a del Pezzo surface of degree $1$ over a number field $k$. Let $S$ be a finite set of places of $k$, and let $v_1,v_2,v_3 \notin S$ be finite places of $k$. Suppose there exists a smooth $\OO_{k,S}$-scheme $\X$ such that the following conditions hold:
\begin{romanenum}
\item $X=\X_k$,
\item $\X_{\FF_{v_1}}$, $\X_{\FF_{v_2}}$, and $\X_{\FF_{v_3}}$ are del Pezzo surface of degree $1$,
\item 
$\phi_{\X_{\FF_{v_1}}}(F_{v_1})$  has order  $7$,
\item 
$\phi_{\X_{\FF_{v_2}}}(F_{v_2})$ has order $6$ and determinant $-1$, and $\phi_{\X_{\FF_{v_2}}}(F_{v_2})^2$ has trace $5$,
\item 
$\phi_{\X_{\FF_{v_3}}}(F_{v_3})$ has order $3$ and trace $-4$.
\end{romanenum}
Then $\phi_X\colon\Gal(\kbar/k) \to \Or(\KXperp)$ is surjective. 
\end{prop}
\begin{proof}
By the assumptions of the proposition and Lemma \ref{L:GaloisLift}, there are $\sigma_1,\sigma_2,\sigma_3 \in \Gal(\kbar/k)$ such that $\phi_X(\sigma_1)$ has order $7$, $\phi_X(\sigma_2)$ has order $6$ and determinant $-1$, $\phi_X(\sigma_2)^2$ has order $3$ and trace $5$, and $\phi_X(\sigma_3)$ has order $3$ and trace $-4$.

Recall that $\KXperp$ is isomorphic to the $E_8$ lattice, and $\Or(\KXperp)\cong W(E_8)$.   Thus applying Proposition~\ref{P:Weyl group criterion}, we find that $\phi_X(\Gal(\kbar/k))=\Or(\KXperp)$.
\end{proof}

\begin{remark}
Let $X$ be a del Pezzo surface of degree $1$ over a number field $k$.  Using the Chebotarev density theorem, it is easy to see that if $\phi_X$ is surjective, then there is a model $\X$ and places $v_1,v_2,v_3$ satisfying the conditions in Proposition~\ref{P:full Galois criterion}.
\end{remark}

Let $W^+(E_8)$ be the subgroup of $W(E_8)$ consisting of the elements with positive determinant.  We have an exact sequence
\begin{equation} \label{E:det SES}
1 \to W^+(E_8) \to W(E_8) \overset{\det}{\to} \{\pm 1\} \to 1.
\end{equation}
Since $-I$ is an element of $W(E_8)$, there is an exact sequence
\begin{equation} \label{E:Schur SES}
1 \to \{\pm I \} \to W^+(E_8) \overset{\varphi}{\to} G \to 1,
\end{equation}
where $G:=W^+(E_8)/\{\pm I\}$ and $\varphi$ is the quotient map.
In \cite{Bourbaki}*{Ch.~VI \S4 Ex.~1}, it is sketched out that $G$ is isomorphic to a certain simple nonabelian group $O^+_8(2)$.  We will use the \emph{Atlas of Finite Groups} \cite{Atlas}*{page~85}, which we will henceforth refer to simply as \emph{the Atlas}, as a source of information concerning the group $G\cong O^+_8(2)$.  In the notation of the Atlas, $W(E_8)$ is isomorphic to $2.O^+_8(2).2$.  

\begin{lemma} \label{L:conjugacy classes of W+}
Given a conjugacy class $C$ of $W^+(E_8)$ of order 3, $\varphi(C)$ is a conjugacy class of $G$ of order 3; this induces a bijection between the conjugacy class of order 3 of $W^+(E_8)$ with those of $G$.

The group $W^+(E_8)$ and $G$ both have exactly five conjugacy classes of order 3.  In both cases these five conjugacy classes have cardinalities $2240$, $2240$, $2240$, $89600$, and $268800$.
\end{lemma}
\begin{proof}
Let $C$ be a conjugacy class of $W^+(E_8)$ of order $3$.  The homomorphism $\varphi$ is surjective, so $\varphi(C)$ is indeed a conjugacy class of $G$.
Since $\ker \phi =\{ \pm I\}$, $\varphi(C)$ must also have order 3.

Let $\Cc$ be a conjugacy class of $G$ of order $3$.  Choose $g \in W^+(E_8)$ such that $\varphi(g)\in \Cc$. Then either $g$ or $(-I)g$ has order $3$ (the other element having order $6$).  
It is now clear that there are two conjugacy classes $C$ of $W^+(E_8)$ such that $\varphi(C)=\Cc$, one of order $3$ and one of order $6$.  The correspondence stated in the lemma is now apparent.

For any conjugacy class $C$ of $W^+(E_8)$ of order $3$, $|C|=|\varphi(C)|$.  The last statement of the lemma on the number and size of conjugacy classes of $G$ of order $3$ can then be read off the Atlas (the Atlas gives the cardinality of the centralizer of any element of a conjugacy class of $G$, from this one can calculate the cardinality of the conjugacy class itself).
\end{proof}

\begin{lemma} \label{L:conjugacy classes of W}
The group $W(E_8)$ has exactly four conjugacy classes of order 3.  These conjugacy classes $\{C_i\}_{i=1,\ldots, 4}$ can be numbered so that 
\begin{align*}
|C_1|&= 2240 , \quad \tr(C_1)= 5   & |C_2|&= 4480, \quad \tr(C_2)= -4 \\
|C_3|&= 89600, \quad \tr(C_3)= -1  & |C_4|&= 268800,\quad \tr(C_4)= 2. 
\end{align*}
\end{lemma}
\begin{proof}
A description of the conjugacy classes of $W(E_8)$ in terms of ``admissible diagrams'' can be found in \cite{MR0318337}.
A conjugacy class of $W(E_8)$ of order 3
must have one of the following four characteristic polynomials: 
\[
\text{$(x^2+x+1)(x-1)^6$, $(x^2+x+1)^2(x-1)^4$, $(x^2+x+1)^3(x-1)^2$ , $(x^2+x+1)^4$}
\]
(and hence have trace $5$, $2$, $-1$, or $-4$ respectively).  In terms of the conventions of \cite{MR0318337}*{\S6}, these characteristic polynomials correspond to the ``admissible diagrams'' $A_2$, $A_2^2$, $A_2^3$ and $A_2^4$.  The lemma is then a consequence of Table 11 in \cite{MR0318337}.
\end{proof}

\begin{proof}[Proof of Proposition \ref{P:Weyl group criterion}]
Define $\Hh = \varphi(H \cap W^+(E_8))$.  Note that any element in $W(E_8)$ of odd order has determinant $+1$.
By assumption (\ref{I:order 7}) of the proposition, there is an $h\in H$ of order $7$.  The homomorphism $\varphi$ has kernel $\{\pm I\}$, so $\varphi(h)$ is an element of order $7$ in $\Hh$.
  
By (\ref{I:negative determinant}) there is a $w\in H$ such that $\det(w)=-1$, and by (\ref{I:trace 5}) and (\ref{I:trace -4}) there are $h_1,h_2\in H\cap W^+(E_8)$ of order $3$ such that $\tr(h_1)=5$ and $\tr(h_2)=-4$.  
Let $\Cc_1$, $\Cc_2$ and $\Cc_3$ be the conjugacy classes of $W^+(E_8)$ of order $3$ and cardinality 2240 (here we are using Lemma \ref{L:conjugacy classes of W+}).  By Lemma \ref{L:conjugacy classes of W} and cardinality considerations, we have (after possibly renumbering the $\Cc_i$) $C_1=\Cc_1$ and $C_2=\Cc_2 \cup \Cc_3$.  

Since $W^+(E_8)$ is a normal subgroup of $W(E_8)$, the set $w\Cc_2 w^{-1}$ is also a conjugacy class of $W^+(E_8)$ (of trace $-4$).  Since $\Cc_2$ is not a conjugacy class of $W(E_8)$, we deduce that $\Cc_3=w\Cc_2 w^{-1}$.  So $h_1 \in C_1=\Cc_1$, and the set $\{h_2,wh_2w^{-1}\}$ contains elements from both of the conjugacy classes $\Cc_2$ and $\Cc_3$. 
Therefore $\varphi(h_1), \varphi(h_2), \varphi(wh_2w^{-1})\in \Hh$ are representatives of the three conjugacy classes of order $3$ and cardinality $2240$ in $G$.\\

Now consider any maximal proper subgroup $M$ of $G$.  The Atlas gives a description of the {maximal} proper subgroups of $G$.  
 
Suppose that $|M| > 155520$.  By checking the permutation character of $G$ associated with $M$ given in the Atlas, one verifies that $M$ does not contain elements from all three of the conjugacy classes of order $3$ and cardinality $2240$ in $G$.  However, we have just shown that $\Hh$ contains elements from each of these three conjugacy classes, hence $\Hh \not \subseteq M$.

 If $|M| \leq 155520$, then the Atlas shows that $7 \nmid |M|$ and in particular $M$ does not have any elements of order $7$.  Since $\Hh$ has an element of order $7$,  $\Hh \not \subseteq M$.
 
Since $\Hh$ is not contained in any of the proper maximal subgroups of $G$, we must have 
\[\varphi(H\cap W^+(E_8))= \Hh = G.\]

Suppose that $-I\not\in H\cap W^+(E_8)$.  Since $\varphi(H\cap W^+(E_8)) = G$ and $(H\cap W^+(E_8))\cap \ker\varphi = \{1\}$, the map $\varphi\colon W^+(E_8)\to G$ induces an isomorphism $H\cap W^+(E_8) \cong G$.  
Using our description of $W(E_8)$ as a subgroup of $\GL(E_8)$, we get an injective group homomorphism $G \hookrightarrow W(E_8)\subseteq \GL(E_8)$.  Hence there exists an \emph{injective} complex representation of $G$,
$$\rho\colon G \hookrightarrow \GL(E_8\otimes_\ZZ \CC)\cong \GL_8(\CC).$$
The character table of $G$ in the Atlas, shows that $G$ has no non-trivial irreducible representations of degree less than $28$, and so any eight dimensional complex representation of $G$ is trivial.  This contradicts the injectivity of $\rho$, and we conclude that $-I$ is an element of $H\cap W^+(E_8)$.

Since $-I \in H \cap W^+(E_8)$ and $\varphi(H \cap W^+(E_8))=G$, the exact sequence (\ref{E:Schur SES}) shows that 
$H\cap W^+(E_8) = W^+(E_8).$
Finally, since $H \supseteq W^+(E_8)$ and $\det(H)=\{\pm 1\}$ (by assumption (\ref{I:negative determinant})), the exact sequence (\ref{E:det SES}) shows that
$H=W(E_8)$.
\end{proof}

\section{Blowup and anticanonical models}
\label{S:dictionary}

\subsection{Equations for anticanonical models} \label{SS:finding the sextic}
\label{SS:equation}
Let $X$ be a del Pezzo surface of degree $1$ over a field $k$.   We now describe how to find a polynomial $f$ as in Proposition~\ref{P:sextic}; details can be found in~\cite[pp.1199--1201]{Cragnolini}.   

Fix a graded $k$-algebra $R=\bigoplus_{m\geq 0} R_m$ that is isomorphic to $R(X,\minusK)$.  In \S\ref{SS:blow-up model}, our ring $R$ will be expressed in terms of the blow-up description of $X$.  By \cite[Corollary III.3.5.2]{Kollar1996}, for each integer $m>0$ 
\begin{equation*}
\dim_k (R_m )= h^0(X,\omega_X^{-m})= \frac{m(m+1)}{2}+1.
\end{equation*} 

\begin{enumerate}
\item Choose a basis $\{x,y\}$ for the $k$-vector space $R_1$. 
\item The elements $x^2$, $xy$, $y^2$ of $R_2$ 
are linearly independent.  Since $\dim_k (R_2)=4$, 
we may choose an element $z$ such that $\{x^2,xy,y^2,z\}$ is a basis of $R_2$.
\item The elements $x^3, x^2y, xy^2, y^3, xz, yz$ of $R_3$ 
are linearly independent.  Since $\dim_k(R_3)=7$, 
we may choose an element $w$ such that $\{x^3,x^2y,xy^2,y^3,xz,yz,w\}$ is a basis of $R_3$.
\item Since $\dim_k(R_6)=22,$
the $23$ elements
\begin{align*}
\{ & x^6,x^5y,x^4y^2,x^3y^3,x^2y^4,xy^5,y^6,x^4z,x^3yz,x^2y^2z,xy^3z, \\
& y^4z,x^2z^2,xyz^2,y^2z^2,z^3,x^3w,x^2yw,xy^2w,y^3w,xzw,yzw,w^2\}
\end{align*}
must be linearly dependent over $k$. Let $f(x,y,z,w) = 0$ be a nonzero linear relation among these elements.  
\item
Viewing $f$ as a sextic polynomial in weighted variables $x,y,z,w$, we have an isomorphism of graded $k$-algebras, $R(X,\minusK)\cong R \cong k[x,y,z,w]/(f)$.
\end{enumerate}

\begin{remark} \label{R:bases}
If $k$ is a field of characteristic not equal to $2$ or $3$, then in step (5) above, we may complete the square with respect to the variable $w$ and the cube with respect to the variable $z$ to obtain an equation involving only the monomials
\[
\{x^6,x^5y,x^4y^2,x^3y^3,x^2y^4,xy^5,y^6,x^4z,x^3yz,x^2y^2z,xy^3z,y^4z,z^3,w^2\}.
\]
\end{remark}

\subsection{The blow-up model}
\label{SS:blow-up model}

\begin{prop} \label{P:dictionary}
Let $k$ be a perfect field.  Fix a $\Gal(\kbar/k)$-stable set $S := \{P_1,\dots,P_8\}$ of eight distinct closed points in $\PP^2_{\kbar} = \Proj(\kbar[x_0,x_1,x_2])$ that are in general position.  Let $\Ii\subseteq \sheafO_{\PP^2_k}$ be the coherent ideal sheaf associated to the closed subset $S$ of $\PP^2_k$ with its reduced-induced subscheme structure.  Let $X$ be the del Pezzo surface of degree $1$ over $k$ obtained by blowing up $\PP^2_k$ along $\Ii$. 
\begin{romanenum}
\item For each $\sigma \in \Gal(\kbar/k)$, the order of $\phi_X(\sigma)$ is equal to the order of the action of $\sigma$ on the set $S$, the trace of $\phi_X(\sigma)$ is equal to the number of $P_i$ fixed by $\sigma$, and the determinant of $\phi_X(\sigma)$ is equal to the sign of the permutation of the $P_i$.
\item There is an isomorphism of graded $k$-algebras,
$$R(X,\minusK) \cong \bigoplus_{m \geq 0} H^0(\PP^2_k, \Ii^m(3m)).$$
The vector space $H^0(\PP^2_k, \Ii^{m}(3m) )$ is the set of homogenous degree $3m$ polynomials in $k[x_0,x_1,x_2]$ that have $m$-fold vanishing at each $P_i$.
\end{romanenum}
\end{prop}
\begin{proof}
First recall that $\Pic(\Xbar)=\KXperp\oplus \ZZ\cdot K_X$, and $\Gal(\kbar/k)$ acts via $\phi_X$ on $\KXperp$ and fixes $K_X$. 

We now give a description of the Galois action in terms of the blow-up data.  Let $\pi \colon X \to \PP^2_{k}$ be the blow-up of $\PP^2_k$ along $\Ii$, and let $e_i$ be the exceptional divisor of $X_{\kbar}$ corresponding to $P_i$.  The action of $\Gal(\kbar/k)$ permutes the exceptional divisors $\{e_i\}$ in the same way that it permutes the points $\{P_i\}$.  In particular, the divisor $e_1+\cdots + e_8$ is defined over $k$. 
Recall that the exceptional divisors $e_1,\dots,e_8$, together with the strict transform $\ell$ of a line in $\PP^2_{\kbar}$ (not passing through any of the $8$ points) form a basis of $\Pic(\Xbar) \cong \ZZ^9$.  The class of ${-K_X} = 3\ell - (e_1+\dots +e_8)$ and $e_1+\dots +e_8$ in $\Pic(\Xbar)$ are Galois stable, so the class of $\ell$ is also stable under $\Gal(\kbar/k)$.   
Therefore $\Pic(\Xbar)=\ZZ e_1\oplus\dots \oplus \ZZ e_8 \oplus \ZZ \ell$, where $\Gal(\kbar/k)$ fixes the class of $\ell$ and permutes the $e_i$ the same way it permutes the $P_i$.  

Part (i) follows by comparing the above descriptions of the action of $\Gal(\kbar/k)$ on $\Pic(\Xbar)$.\\
 
Let $\tilde{\Ii}\subseteq \OO_X$ be the ideal sheaf of $e_1\cup\dots\cup e_8$ with the reduced-induced subscheme structure.   We claim that the invertible sheaf
\[
\tilde{\Ii}\otimes_{\OO_X} \pi^*(\OO_{\PP_k^2}(3))
\]
is isomorphic to the anticanonical sheaf of $X$.   To prove this it suffices to work over an algebraic closure $\kbar$; in which case the invertible sheaves $\pi^*(\OO_{\PP_k^2}(3))$ and $\tilde{\Ii}$ correspond to the divisor classes of $3\ell$ and $-(e_1+\dots+e_8)$ in $\Pic(\Xbar)$, respectively.   The claim is then immediate since $-K_X = 3\ell - (e_1+\dots + e_8) \in \Pic(\Xbar)$.

We thus have the following isomorphisms of graded $k$-algebras (we shall write $\PP^2$ instead of $\PP^2_k$ to avoid clutter), 
\begin{align*}
R(X,\minusK) &\cong \bigoplus_{m\geq 0} H^0(X,(\tilde{\Ii}\otimes_{\OO_X} \pi^*(\OO_{\PP^2}(3)))^{\otimes m}) \\
& \cong \bigoplus_{m\geq 0} H^0(X,\tilde{\Ii}^m\otimes_{\OO_X} \pi^*(\OO_{\PP^2}(3m))) \\
& = \bigoplus_{m\geq 0} H^0(\PP^2,\pi_*(\tilde{\Ii}^m\otimes_{\OO_X} \pi^*(\OO_{\PP^2}(3m)))) \\
& \cong \bigoplus_{m\geq 0} H^0(\PP^2,\pi_*(\tilde{\Ii}^m)(3m)),
\end{align*}
where the last isomorphism follows from the projection formula \cite{Hartshorne}*{Exercise II.5.1}.

The morphism $\OO_{\PP^2} \to \pi_*\OO_X$ coming from $\pi$, is an isomorphism of $\OO_{\PP^2}$-modules (this can be checked on stalks, using the fact that $\pi$ gives an isomorphism between $X - \pi^{-1}(S)$ and $\PP^2 - S$, and that a regular function on $e_i$ must be constant).   We will identify $\OO_{\PP^2}$ and $\pi_*\OO_X$, and in particular, we may view $\pi_*(\tilde{\Ii})$ as an ideal of $\OO_{\PP^2}$.  

The ideal sheaves $\pi_*(\tilde{\Ii})$ and $\Ii$ both correspond to closed subschemes of $\PP^2$ with support $\{P_1,\dots,P_8\}$.  Since $\Ii$ corresponds to the reduced-induced subscheme structure on $\{P_1,\dots,P_8\}$, we find that $\pi_*(\tilde{\Ii}) \subseteq \Ii.$

The ideal sheaf $\pi^*(\Ii)$ has support $e_1\cup \dots \cup e_8.$  Since $\tilde{\Ii}$ corresponds to the reduced-induced subscheme structure on $e_1\cup \dots \cup e_8$, we find that $\pi^*(\Ii) \subseteq \tilde{\Ii}$ and hence
\[
\Ii= \pi_*\pi^*(\Ii) \subseteq \pi_*(\tilde{\Ii})
\]
(the equality follows from the projection formula). Therefore $\pi_*(\tilde{\Ii}) = \Ii$.   For any $m\geq 1$, 
we have $\pi_*(\tilde{\Ii}^m) =\pi_*(\tilde{\Ii})^m$ (these sheaves have support $\{P_1,\dots,P_8\}$ so it suffices to check on the stalks at each $P_i$).  Thus $\pi_*(\tilde{\Ii}^m)= \pi_*(\tilde{\Ii})^m = \Ii^m$ for each $m\geq 1$, and hence
\[
R(X,\minusK) \cong \bigoplus_{m\geq 0} H^0(\PP^2,{\Ii}^m(3m)).   \qedhere
\]
\end{proof}

The graded ring in Proposition~\ref{P:dictionary}(ii) is amenable to computation, and in particular we may implement the procedure from \S\ref{SS:finding the sextic}.  Thus, given a Galois stable set of $8$ points in general position, we have a method for finding the corresponding del Pezzo surface of degree $1$ as a sextic in weighted projective space.  We will give two examples of this. The first will be of a surface $X$ over $\FF_3$ such that $\phi_{X}(F_3)$ has order $7$. The second will be a surface $X$ over $\FF_5$ such that $\phi_X(F_5)$ has order $6$ and determinant $-1$, and such that $\phi_X(F_5)^2$ has trace $5$. We will use these surfaces to construct a del Pezzo surface of degree $1$ over $\QQ$ satisfying parts (iii) and (iv) of Proposition~\ref{P:full Galois criterion}. The calculations amount simply to linear algebra over finite fields, and are easily implemented on a computer.  We have provided enough details so that the careful reader may verify all our claims.
 
\subsection{An example over $\FF_3$}
\label{Ex:F3 example}

Let $k = \FF_3$, and let $\alpha$ be a root of the irreducible polynomial $t^7 + 2t^2 + 1\in \FF_3[t]$.  Let $F_3\colon \FFbar_3 \to \FFbar_3$ be the Frobenius map $x\mapsto x^3$, and define the following set of eight points in $\PP^2_{\overline{\FF}_3}$, which are in general position:
$$
S := \{[1,0,0]\} \cup \{[1:F_3^i(\alpha):F_3^i(\alpha^4)] : 0\leq i \leq 6 \}.
$$
Let $\Ii$ be the coherent sheaf of ideals in $\PP^2_{\FF_3}$ corresponding to the reduced-induced structure of $S$. Denote by $\pi\colon X\to \PP^2_{\FF_3}$ the blow-up along $\Ii$. Since $F_3$ fixes $[1:0:0]$ and acts transitively on the other seven points of $S$, it follows from Proposition~\ref{P:dictionary}(i) that $\phi_X(F_3)$ has order $7$.  

We now use Proposition~\ref{P:dictionary}(ii) and the procedure of \S\ref{SS:equation} to calculate a defining equation for $X$.  The polynomials
$$
x = x_0^2x_1 + x_0x_2^2 + 2x_1^3 \quad\text{and}\quad 
y = x_0^2x_2 + 2x_0x_1^2 + 2x_1x_2^2
$$
form a basis of $H^0(\PP^2_{\FF_3},\Ii(3))$. Let
$$
z = x_0^6 + x_0^3x_1^3 + 2x_0^2x_1^2x_2^2 + 2x_0^2x_2^4 + x_0x_1^5 + 2x_0x_1^3x_2^2 + x_0x_1^2x_2^3 + 
    x_0x_1x_2^4 + 2x_1^6 + 2x_1^5x_2 + x_1^4x_2^2;
$$
then $\{x^2,xy,y^2,z\}$ is a basis of $H^0(\PP^2_{\FF_3},\Ii^2(6))$. Let
\begin{align*}
w &= 
x_0^9 + 2x_0^7x_1x_2 + x_0^5x_1^2x_2^2 + x_0^5x_2^4 + 2x_0^4x_1^2x_2^3 + 2x_0^4x_1x_2^4 + x_0^3x_1^6 
    + x_0^3x_1^4x_2^2 + x_0^3x_2^6 \\
    &\quad + 2x_0^2x_1^6x_2 + 2x_0^2x_1^4x_2^3 + x_0x_1^7x_2 + x_0x_1^6x_2^2 
    + x_0x_1^5x_2^3 + 2x_0x_1^2x_2^6 + x_1^5x_2^4 + 2x_1^4x_2^5 + 2x_1^3x_2^6;
\end{align*}
then $\{x^3,x^2y,xy^2,y^3,xz,yz,w\}$ is a basis of $H^0(\PP^2_{\FF_3},\Ii^3(9))$. Elementary linear algebra now yields the relation
\[
2x^6 + 2x^3y^3 + xy^5 + y^6 + 2x^3yz + x^2y^2z + y^4z + 2xyz^2 + 2z^3 + 2x^3w  + y^3w + 
    w^2 = 0
\]
in $H^0(\PP^2_{\FF_3},\Ii^6(18))$. Completing the square with respect to $w$, we obtain the equation
\begin{equation*}
x^6 + x^3y^3 + xy^5 + 2x^3yz + x^2y^2z + y^4z + 2xyz^2 + 2z^3 + w^2 = 0.
\end{equation*}
Under the transformation $z\mapsto -z$ the coefficient of $z^3$ becomes $1$ and yields
\begin{equation}
\label{E:order7element}
x^6 + x^3y^3 + xy^5 + x^3yz + 2x^2y^2z + 2y^4z + 2xyz^2 + z^3 + w^2 = 0.
\end{equation}
This gives a model for $X$ as a smooth sextic hypersurface in $\PP_{\FF_3}(1,1,2,3)$. 

\subsection{An example over $\FF_5$}
\label{Ex:F5 example}

Let $k = \FF_5$, and let $\alpha$ be a root of the irreducible polynomial $t^6 + t^4 + 4t^3 + t^2 + 2 \in \FF_5[t]$.  Define $\beta = \alpha^{5^4 + 5^2 + 1}$ and $\gamma = \alpha^{5^3 + 1}$, so $\FF_5(\beta)$ and $\FF_5(\gamma)$ are degree $2$ and $3$ extensions of $\FF_5$, respectively.  Let $F_5\colon \FFbar_5 \to \FFbar_5$ be the Frobenius map $x\mapsto x^5$, and define the following set of eight points in $\PP^2_{\overline{\FF}_5}$, which are in general position:
\begin{align*}
S &:= \{[1,0,0],[3:2:4],[4:2:1],[1:\beta:\beta^3],[1:F_5(\beta):F_5(\beta^3)], \\
&\qquad [1:\gamma:\gamma^4],[1:F_5(\gamma):F_5(\gamma^4)],[1:F_5^2(\gamma):F_5^2(\gamma^4)]\}.
\end{align*}
Let $\Ii$ be the coherent sheaf of ideals in $\PP^2_{\FF_5}$ corresponding to the reduced-induced structure of $S$. Denote by $\pi\colon X\to \PP^2_{\FF_5}$ the blow-up along $\Ii$. Since $F_5$ acts as an order $6$ odd permutation on $S$, it follows from Proposition~\ref{P:dictionary}(i) that $\phi_X(F_5)$ has order $6$ and determinant $-1$.  The automorphism $F_5^2$ fixes exactly 5 elements of $S$, that is the three $\FF_5$-rational points and the order $2$ orbit of $[1:\beta:\beta^3]$; therefore $\phi_X(F_5)^2$ has trace $5$.  

We now use Proposition~\ref{P:dictionary}(ii) and the procedure of \S\ref{SS:equation} to calculate a defining equation for $X$.  The polynomials
\begin{align*}
x &= x_0^3 + 4x_0x_1^2 + 2x_0x_1x_2 + x_0x_2^2 + x_1^2x_2 + 4x_1x_2^2 \quad\text{and}\quad  \\
y &= x_0^2x_1 + 3x_0^2x_2 + 3x_0x_1^2 + x_0x_1x_2 + 3x_0x_2^2 + 4x_1^2x_2 + 3x_1x_2^2 \\
\end{align*}
form a basis of $H^0(\PP^2_{\FF_5},\Ii(3))$. Let
\begin{align*}
z &= x_0^5x_1 + 2x_0^4x_2^2 + 4x_0^3x_1^3 + 2x_0^3x_1^2x_2 + x_0^3x_1x_2^2 + 4x_0^3x_2^3 + 
    3x_0^2x_1^4 + x_0^2x_1^3x_2 + 4x_0^2x_1^2x_2^2 \\
    &\quad + 3x_0^2x_1x_2^3 + 3x_0^2x_2^4 + 
    x_0x_1^4x_2 + 4x_0x_1^3x_2^2 + 4x_0x_1^2x_2^3 + 3x_0x_2^5 + 4x_1^4x_2^2 + 4x_1^3x_2^3 \\
    &\quad + x_1^2x_2^4 + 2x_1x_2^5 + 3x_2^6;
\end{align*}
hen $\{x^2,xy,y^2,z\}$ is a basis of $H^0(\PP^2_{\FF_5},\Ii^2(6))$. Let
\begin{align*}
w &= 
x_0^9 + 2x_0^6x_1^2x_2 + 2x_0^6x_1x_2^2 + x_0^6x_2^3 + x_0^5x_1^3x_2 + 3x_0^5x_1^2x_2^2 + 
    4x_0^5x_1x_2^3 + 3x_0^5x_2^4 + 4x_0^4x_1^5 \\
    &\quad + 3x_0^4x_1^3x_2^2 + 2x_0^4x_2^5 + 
    3x_0^3x_1^5x_2 + x_0^3x_1^4x_2^2 + 3x_0^3x_1^3x_2^3 + 3x_0^3x_1^2x_2^4 + 4x_0^3x_1x_2^5
    + 3x_0^3x_2^6 \\
    &\quad + 2x_0^2x_1^5x_2^2 + x_0^2x_1^4x_2^3 + x_0^2x_1^3x_2^4 + 4x_0^2x_1^2x_2^5 +
    4x_0^2x_1x_2^6 + 2x_0^2x_2^7 + 2x_0x_1^6x_2^2 + 4x_0x_1^5x_2^3 \\
    &\quad + x_0x_1^4x_2^4 + 
    3x_0x_1^3x_2^5 + 4x_0x_1^2x_2^6 + 2x_0x_1x_2^7 + x_1^6x_2^3 + x_1^5x_2^4 + x_1^4x_2^5 + 
    2x_1^2x_2^7 + x_2^9;
    \end{align*}
then $\{x^3,x^2y,xy^2,y^3,xz,yz,w\}$ is a basis of $H^0(\PP^2_{\FF_5},\Ii^3(9))$. Elementary linear algebra now yields the linear relation
\begin{align*}
2x^6 &+ 3x^5y + x^4y^2 + 4x^3y^3 + 4x^2y^4 + 4y^6 + 4x^4z + 2x^3yz + x^2y^2z + 2xy^3z \\
    &+ 3y^4z + 3x^2z^2 + 2y^2z^2 + 2z^3 + 2x^3w + 2x^2yw + 2xy^2w + xzw + w^2 = 0
\end{align*}
in $H^0(\PP^2_{\FF_5},\Ii^6(18))$.
Completing the square with respect to $w$ and the cube with respect to $z$, we obtain the equation
\begin{equation*}
2x^5y + x^4y^2 + 2x^3y^3 + 3x^2y^4 + xy^5 + 2x^4z + x^3yz + 4x^2y^2z + 2xy^3z + 4y^4z
    + 2z^3 + w^2 = 0.
\end{equation*}
Multiplying both sides by $4$, and rescaling by $[x,y,z,w]\mapsto [x,y,z/2,w/2]$ yields
\begin{equation}
\label{E:twobirds}
3x^5y + 4x^4y^2 + 3x^3y^3 + 2x^2y^4 + 4xy^5 + 4x^4z + 2x^3yz + 3x^2y^2z + 4xy^3z + 3y^4z + z^3 + w^2 = 0.
\end{equation}
This gives a model for $X$ as a smooth sextic hypersurface in $\PP_{\FF_5}(1,1,2,3)$. 

\section{The Lefschetz trace formula}

In \S\ref{SS:F7 example}, we will describe a del Pezzo surface $X$ of degree $1$ over $\FF_7$ such that $\phi_X(F_7)$ has order $3$ and trace $-4$.  This gives a surface satisfying part (v) of Proposition~\ref{P:full Galois criterion}.   To prove these properties of $X/\FF_7$ it will suffice, by the Lefschetz trace formula, to compute $|X(\FF_7)|$ and $|X(\FF_{7^3})|$.  The methods used in \S\ref{S:dictionary}  cannot produce this example, since Proposition~\ref{P:dictionary} only gives surfaces $X/\FF_7$ with $\tr(\phi_X(F_7))\geq 0$.  

\subsection{The Lefschetz trace formula}
The following version of the Lefschetz trace formula (specialized to del Pezzo surfaces) is due to Weil.  
\begin{thm} \label{T:Lefschetz}
Let $\FF_q$ be a finite field with $q$ elements, and let $F_q \in \Gal(\FFbar_q/\FF_q)$ be the Frobenius automorphism $x\mapsto x^q$.  Let $X$ be a del Pezzo surface over $\FF_q$ of degree $d\leq 6$.   Then
\[
|X(\FF_q)| = q^2 + q(\tr(\phi_X(F_q))+1)+ 1.
\]
\end{thm}
\begin{remark}
For a proof of Theorem~\ref{T:Lefschetz}, see \cite{Manin}*{\S27}.  Note that $\tr(\phi_X(F_q))+1$ is the trace of the action of $F_q$ on $\Pic(X_{\FFbar_q})$.
\end{remark}

\subsection{Points on the anticanonical model} 
\label{SS:points on anticanonical model}
Fix a finite field $\FF$.  Let $X$ be a del Pezzo surface of degree $1$ defined over $\FF$ which is given explicitly as a smooth sextic hypersurface
\[
w^2= z^3 + F(x,y)z^2 + G(x,y)z + H(x,y)
\]
in $\PP_{\FF}(1,1,2,3)$.  
Now consider the morphism $\varphi\colon X-\{[0:0:1:1]\} \to \PP^1_{\FF}$, $[x:y:z:w] \mapsto [x:y]$ of Lemma~\ref{L:anticanonical map}.  
Take any point $P=[a:b]\in\PP^1(\FF)$ with $(a,b)\in\FF^2-\{(0,0)\}$.  The fiber of $\varphi$ above $P$ is isomorphic to the affine curve
\[
C_{P}\colon W^2 = Z^3 + F(a,b)Z^3 + G(a,b)Z^2 +H(a,b)
\]
in $\AA^2_{\FF}$.  We thus have
\begin{equation} \label{E:how to count}
|X(\FF)| = \sum_{P\in \PP^1(\FF)} |C_{P}(\FF)| + 1.
\end{equation}

\subsection{Example over $\FF_7$}
\label{SS:F7 example}
\begin{lemma}
\label{L:F7 example}
Let $X$ be the closed subscheme of $\PP_{\FF_7}(1,1,2,3)$ defined by the sextic
\[
w^2=z^3 + 2x^6 + 2y^6.
\]
Then $X$ is a del Pezzo surface of degree $1$ over $\FF_7$, $\phi_X(F_7)$ has order $3$, and $\tr (\phi_X(F_7)) = -4$.
\end{lemma}
\begin{proof}
The scheme $X$ is defined by a smooth sextic, and hence is a del Pezzo surface of degree $1$ by Proposition~\ref{P:sextic}.  Consider an element $g\in \Or(\KXperp)$.  Since $\KXperp\cong \ZZ^8$ and $g$ has finite order, we find that $\tr(g)\leq 8$, with equality holding if and only if $g=I$.  By Theorem~\ref{T:Lefschetz},
\[
|X(\FF_7)| = 7^2 + 7  (\tr(\phi_X(F_7))+1) + 1 \text{\quad and \quad} |X(\FF_{7^3})| = 7^6 + 7^3 (\tr(\phi_X(F_7)^3)+1) + 1,
\]
and thus the lemma is equivalent to showing that
\[
|X(\FF_7)| = 7^2 + 7\cdot(-3) + 1 = 29\text{\quad and \quad} |X(\FF_{7^3})| = 7^6 + 7^3 \cdot 9 + 1 = 120737.
\]
Let $\FF$ be an extension of $\FF_7$.  For $(a,b)\in \FF^2-\{(0,0)\}$, define the affine curve
\[
C_{[a,b]}\colon W^2 = Z^3 + 2a^6+2b^6
\]
in $\AA^2_{\FF}$.  From (\ref{E:how to count}),
\begin{align} \label{E:how to count, fibered}
|X(\FF)| & = |C_{[1,0]}(\FF)| + \sum_{a\in \FF} |C_{[a,1]}(\FF)| + 1.
\end{align}
For $\FF=\FF_7$, we have $|C_{[1,0]}(\FF_7)|=|C_{[0,1]}(\FF_7)|=8$, and $|C_{[a,1]}(\FF_7)|= 2$ for all $a\in\FF_7^\times$; so
$|X(\FF_{7})| = 2\cdot 8 + 6\cdot 2 + 1 = 29.$

Now let $\FF=\FF_{7^3}$.  We will use (\ref{E:how to count, fibered}) to compute $|X(\FF)|$, but it is useful to note that $|C_{[a,b]}(\FF)|$ depends only the class of $2a^6+2b^6$ in $\FF^\times/(\FF^\times)^6 \cup \{0\}$.  We have a bijection $r\colon \FF^\times/(\FF^\times)^6 \cup\{0\} \to \FF_7$, $a\mapsto a^{(|\FF|-1)/6}$.  For $i\in\FF_7$, let $N_i(\FF)$ be the the number of $\FF$-points of the affine curve $W^2=Z^3 + \alpha$, where $\alpha$ is any element of $\FF$ with $r(\alpha)=i$. It follows that
\[
|X(\FF)|  = |C_{[1,0]}(\FF)| + \sum_{i=0}^6 |\{a \in \FF: r(2a^6+2)=i\}|\cdot N_i(\FF) + 1.
\]
The right hand side is readily computed, and we find that:
\begin{align*}
|X(\FF_{7^3})| & =  323 + 0 \cdot 343 + 43\cdot 323 + 72\cdot 380 + 72\cdot 360 + 36\cdot 326 + 36\cdot 306 + 84\cdot 363 + 1 \\
&= 120737.  \qedhere
\end{align*}
\end{proof}

\begin{remark}
We can also verify our previous examples of del Pezzo surfaces over finite fields using this method.  For example, consider the surface $X/\FF_3$ defined by (\ref{E:order7element}).  To show that $\phi_X(F_3)$ has order $7$, it suffices to check that
\[
|X(\FF_3)| \neq 3^2 + 3\cdot 9 +1 \text{\quad and \quad} |X(\FF_{3^7})| = 3^{14} + 3^7\cdot 9 + 1.
\]
Now consider the surface $X/\FF_5$ defined by (\ref{E:twobirds}).  One verifies that
\[
 |X(\FF_{5^2})| = 5^4 + 5^2\cdot 6 + 1, \quad |X(\FF_{5^3})| = 5^6 + 5^3\cdot 7 +1 , \text{\quad and \quad} |X(\FF_{5^6})| = 5^{12} + 5^6\cdot 9 + 1.
\]
It is then apparent that $\phi_X(F_5)$ has order $6$ and $\phi_X(F_5)^2$ has trace $5$.  Since $\phi_X(F_5)^3$ has order $2$ and $\tr(\phi_X(F_5)^3)=6$, we deduce that $\phi_X(F_5)^3$ has eigenvalue $+1$ with muliplicity $7$, and $-1$ with multiplicity $1$.  Therefore $\det(\phi_X(F_5))=\det(\phi_X(F_5)^3) = -1$. 
\end{remark}

\section{Proof of Theorem~\ref{T:main theorem}}  \label{S:main proof}

\noindent Let $\X = \Proj(\ZZ[x,y,z,w]/(f))$; since 
\[
f \equiv x^6 + x^3y^3 + xy^5 + x^3yz + 2x^2y^2z + 2y^4z + 2xyz^2 +z^3 + w^2 \smallpmod{3},
\] 
we find that $\X_{\FF_3}$ is the del Pezzo surface of degree $1$ from Example~\ref{Ex:F3 example}.  In particular, $\phi_{\X_{\FF_3}}(F_3)$ has order $7$.  Since 
\begin{align*}
f \equiv 3x^5y &+ 4x^4y^2+3x^3y^3+ 2x^2y^4+4xy^5+4x^4z\\ &+2x^3yz+3x^2y^2z + 4xy^3z + 3y^4z +z^3+w^2 \smallpmod{5},
\end{align*}
we find that $\X_{\FF_5}$ is the del Pezzo surface of degree $1$ from Example~\ref{Ex:F5 example}.  In particular, $\phi_{\X_{\FF_5}}(F_5)$ has order $6$ and determinant $-1$, and $\phi_{\X_{\FF_5}}(F_5)^2$ has trace $5$. 
Since 
\[
f\equiv 5x^6 + 5y^6 + z^3 + w^2 \smallpmod{7},
\] 
we find that $\X_{\FF_7}$ is isomorphic to the del Pezzo surface of degree $1$ from Lemma~\ref{L:F7 example}; thus $\phi_{\X_{\FF_7}}(F_7)$ has order $3$ and trace $-4$.

Let $S$ be a finite set of primes such that $\X' := \X_{\ZZ[S^{-1}]}$ is smooth over $\Spec\ZZ[S^{-1}]$, where $\ZZ[S^{-1}]$ is the ring of $S$-units in $\QQ$.  Since $\X$ has smooth fibers at $3,5$ and $7$, we may assume that $S$ is chosen such that $3, 5, 7 \notin S$.  Note that $X = \X'_{\QQ}$ is smooth, so it is a del Pezzo surface of degree $1$ by Proposition~\ref{P:sextic}.  By Proposition~\ref{P:full Galois criterion}, we deduce that $\phi_X\colon\Gal(\Qbar/\QQ) \to \Or(\KXperp)$ is surjective.\\

To prove the final statement of the theorem, it suffices to show that for any number field $k\subseteq \Qbar$, we may find an $f$ such that $L_X\cap k=\QQ$.  
Let $k_1,\dots, k_m$ be all the subfields of $k$ except for $\QQ$.  By the Chebotarev density theorem, there are distinct rational primes $p_1,\dots,p_m$ greater than $7$ such that $p_i$ does not split completely in $k_i$.

For each $i$, choose 8 points of $\PP^2(\FF_{p_i})$ that are in general position; blowing them up gives a del Pezzo surface $X_i$ of degree $1$ defined over $\FF_{p_i}$.  By Proposition~\ref{P:dictionary}(i), $\phi_{X_i} = 1$.

There is a sextic polynomial $f_i \in \FF_{p_i}[x,y,z,w]$ such that $X_i$ is isomorphic to the hypersurface $f_i=0$ in $\PP_{\FF_{p_i}}(1,1,2,3)$ (the polynomial can be calculated using the results of \S\ref{S:dictionary}).
Now let $f\in \ZZ[x,y,z,w]$ be a sextic which satisfies (\ref{E:main equation}), and for each $i$ satisfies $f\equiv f_i \smallpmod{p_i}$.  There is a finite set $S$ of rational primes such that the scheme $\X := \Proj(\ZZ[S^{-1}][x,y,z,w]/(f))$ is smooth over $\Spec\ZZ[S^{-1}]$, and $S$ can be chosen so that $3, 5, 7, p_1,\dots, p_m \notin S$. We have already proven that $X := \X_\QQ$ is a del Pezzo surface of degree $1$ with surjective homomorphism $\phi_X$.  

If $L_X\cap k\neq \QQ$, then $k_i=L_X\cap k$ for some $i$. Let $\Fr_{p_i}\in\Gal(\Qbar/\QQ)$ be any Frobenius automorphism over $p_i$.  Using Lemma~\ref{L:specialization} and $\phi_{\X_{p_i}}(F_{p_i}) = \phi_{X_i}(F_{p_i}) =1$, one proves that $\phi_X(\Fr_{p_i})=1$.  Therefore $p_i$ splits completely in $L_X$ (and hence also in $L_X\cap k=k_i$).  This contradicts the assumption that $p_i$ does not split completely in $k_i$.   We conclude that $L_X\cap k=\QQ$.

\section{Elliptic curves}\label{S:Elliptic curve proofs}
\subsection{Rational elliptic surfaces}
\label{S:MWLs}

We summarize, making certain simplifying assumptions, some basic facts about Mordell-Weil lattices of rational elliptic surfaces. A full account of the theory can be found in~\cite{Shioda2}.

Let $\pi\colon\calE \to \PP^1_\Qbar$ be an elliptic surface with a fixed section $\OO$. Assume that 
\begin{enumerate}
\item $\calE$ is rational
\item $\pi$ has at least one singular fiber, and no reducible fibers (in Shioda's notation, $R = \emptyset$).
\end{enumerate}

Let $E$ be the generic fiber of $\pi\colon\calE \to \PP^1_\Qbar$, which is an elliptic curve over the function field $\Qbar(t)$.  There is a natural one-to-one correspondence between the $\Qbar(t)$-points of $E$ and the sections of $\pi\colon\calE \to \PP^1_\Qbar$. The image of the section corresponding to a point $P \in E(\Qbar(t))$ will be denoted by $(P)$; it is a divisor of the surface $\calE$. 

To each point $P\in E(\Qbar(t))$ we associate a fibral divisor\footnote{An irreducible divisor $\Gamma$ of $\calE$ is \defi{fibral} if $\pi|_\Gamma\colon \Gamma \to \PP^1_{\Qbar}$ is a constant map.} $\Phi_P \in \Div(\calE)\otimes_\ZZ \QQ$ such that for all fibral $F \in \Div(\calE)$
\[
((P) - (\OO) + \Phi_P,F) = 0
\]
(such a divisor always exists, see \cite[III.8.3]{Silverman2}).  Let $\NS(\calE)$ be the N\'eron-Severi group of $\calE$; in our case, this group is finitely generated and torsion free (\cite[Theorem 1.2]{Shioda2}).  The map
\begin{align*}
\phi \colon E(\Qbar(t)) &\to \NS(\calE)\otimes_\ZZ \QQ \\
P &\mapsto (P) - (\OO) + \Phi_P
\end{align*}
is a group homomorphism with kernel $E(\Qbar(t))_{\tors}$~\cite[8.2]{Shioda2}. We define a pairing on $E(\Qbar(t))$ using the intersection pairing of $\calE$ as follows:
\begin{align*}
\langle\,\cdot\, ,\cdot\, \rangle \colon E(\Qbar(t)) \times E(\Qbar(t))\to \QQ,\quad 
\langle P,Q\rangle = -(\phi(P),\phi(Q)).
\end{align*}
This pairing is symmetric, bilinear, and coincides with the canonical height pairing on $E(\Qbar(t))$~\cite[III.9.3]{Silverman2}.

Let $T$ be the subgroup of $\NS(\calE)$ generated by $(\OO)$ and all the fibers of $\pi\colon\calE\to \PP^1_{\Qbar}$. By~\cite[Theorem 1.3]{Shioda2}, we have a group isomorphism
\[
\beta\colon E(\Qbar(t)) \to \NS(\calE)/T,\quad P\mapsto (P) \bmod T.
\]
Let $T' = (T\otimes \QQ)\cap \NS(\calE)$. We have an isomorphism $\beta\colon E(\Qbar(t))_{\tors} \overset{\sim}{\to} T'/T$ by~\cite[Corollary 5.3]{Shioda2}. Hence there is an isomorphism of lattices
\begin{equation}
\label{E:LatticeIso}
\beta\colon E(\Qbar(t))/E(\Qbar(t))_{\tors} \overset{\sim}{\to} \NS(\calE)/T'.
\end{equation}
\begin{remark}
The isomorphism~(\ref{E:LatticeIso}) holds without the hypothesis that $\pi\colon\calE \to \PP^1_\Qbar$ has no reducible fibers, but the map $\phi$ is harder to define in this case.
\label{R:reducible}
\end{remark}

\subsection{Galois actions}
\label{S:galoisactions}
Let $\pi\colon\calE \to \PP^1_\QQ$ be an elliptic surface with a fixed section $\OO$, such that the elliptic surface $\bar{\pi}\colon\calE_\Qbar\to \PP^1_\Qbar$, obtained by base extension, satisfies the hypotheses (1) and (2) of \S\ref{S:MWLs}. Then the isomorphism of lattices~(\ref{E:LatticeIso}) respects the $\Gal(\Qbar/\QQ)$-actions~\cite[proof of 8.13]{Shioda2}. Hence $E(\Qbar(t))/E(\Qbar(t))_{\tors}$ and $\NS(\calE)/T'$ are also isomorphic as $\Gal(\Qbar/\QQ)$-modules.

\subsection{Elliptic surfaces associated to del Pezzo surfaces of degree $1$}
\label{S:elliptic surfaces on dP1s}

Let $X$ be a del Pezzo surface of degree $1$ over $\QQ$.  We now describe how, given $X$, one can obtain a rational elliptic surface. The linear system $|{-K}_X|$ gives rise to a rational map $f\colon X \dashrightarrow \PP^1_\QQ$ that is regular everywhere except at the anticanonical point $O$ (cf. Lemma~\ref{L:anticanonical map}).  Blowing up $X$ at $O$, we obtain a surface $\calE$.  Composing the blow-up map with $f$ gives a morphism $\pi\colon \calE \to \PP^1_\QQ$, where almost all of the fibers are non-singular genus $1$ curves.  The morphism $\pi$ induces an isomorphism between the exceptional divisor of $\calE$ corresponding to $O$ and $\PP^1_\QQ$; we thus have a distinguished section $\OO\colon \PP^1_\QQ \to \calE$ of $\pi$.  Therefore, $\pi\colon \calE \to \PP^1_\QQ$ with the section $\OO$ is an elliptic surface.

Concretely, if $X$ is given by a smooth sextic
\[
w^2 = z^3 + F(x,y)z^2 + G(x,y)z + H(x,y)
\]
in $\PP_\QQ(1,1,2,3)$, then the anticanonical point is $O = [0:0:1:1]$.  In this case, $\calE$ is the subscheme of $\PP_{\QQ}(1,1,2,3)\times \PP^1_{\QQ} = \Proj(\QQ[x,y,z,w])\times\Proj(\QQ[u,v])$ cut out by the equations
\begin{equation}
\label{E:blow-up model}
w^2 = z^3 + F(x,y)z^2 + G(x,y)z + H(x,y)\quad \text{and}\quad vx - uy.
\end{equation}
The map $\pi\colon \calE \to \PP_\QQ^1$ is then given by $([x:y:z:w],[u:v]) \mapsto [u:v]$. Note that for points away from the exceptional divisor we have $[u:v] = [x:y]$.

Let $t$ be the rational function $u/v$, thus $x = ty$ on $\calE$. The generic fiber of $\pi$ is the curve  \begin{equation}
\label{E:weighted generic fiber}
E\colon w^2 = z^3 + y^2F(t,1)z^2 + y^4G(t,1)z + y^6H(t,1)
\end{equation}
in $\Proj(\QQ(t)[y,z,w])$. On the \emph{affine} chart $\Spec(\QQ(t)[z/y^2,w/y^3])$ of this weighted ambient space, the curve~(\ref{E:weighted generic fiber}) is isomorphic to the affine curve
\[
(w/y^3)^2 = (z/y^2)^3 + F(t,1)(z/y^2)^2 + G(t,1)(z/y^2) + H(t,1).
\]
Relabelling the variables, we find that the elliptic curve $E/\QQ(t)$ is given by the Weierstrass model
\[
y^2 = x^3 + F(t,1)x^2 + G(t,1)x + H(t,1).
\]

\subsection{Proof of Theorem \ref{T:curvewithsymmetry}}
Let $X$ be the a del Pezzo surface as in Theorem~\ref{T:main theorem}.  
Let $\pi\colon\calE\to\PP^1_\QQ$ be the elliptic surface obtained by blowing up the anticanonical point of $X$ (see \S\ref{S:elliptic surfaces on dP1s}).  The generic fiber of this surface is the elliptic curve $E/\QQ(t)$ in the statement of the theorem. Let $\bar{\pi}\colon\calE_\Qbar \to \PP^1_\Qbar$ be the base extension of $\pi$ by $\Spec\Qbar\to\Spec\QQ$.  The following properties hold:
\begin{enumerate} 
\item
the surface $\calE_\Qbar$ is rational (since it is isomorphic to a blow-up of $\PP_\Qbar^2$ at $9$ points).
\item 
$\bar{\pi}$ has at least one singular fiber (otherwise $\calE$ has constant $j$-invariant~\cite[Ex.~3.35(c)]{Silverman2}). Moreover, the fibers of $\bar{\pi}\colon\calE_\Qbar \to \PP^1_\Qbar$ are irreducible: using the blow-up model~(\ref{E:blow-up model}) of $\calE$, the reader may verify that the fiber above the point $[u:v] \in \PP^1_\Qbar$ is isomorphic to the projectivization of the irreducible curve
\[
y^2 = x^3 + F(u,v)x^2 + G(u,v)x + H(u,v).
\]
\end{enumerate}
We may thus apply \S\ref{S:MWLs} to obtain an isomorphism
\[
E(\Qbar(t))/E(\Qbar(t))_{\tors}  \overset{\sim}{\to} \NS(\calE_\Qbar)/T'.
\]
Furthermore, this isomorphism respects the action of $\Gal(\Qbar/\QQ)$ (see \S\ref{S:galoisactions}). On the other hand, the lattice $\KXperp$ is isomorphic to $\NS(\calE_\Qbar)/T'$ via the composition of maps
\[
\KXperp \to \NS(\calE_\Qbar) \to \NS(\calE_\Qbar)/T',
\]
where the first map is induced by pullback of divisors along the blow-up map, and the second is the natural quotient map.  Therefore we have isomorphisms
\[
\KXperp \cong \NS(\calE_\Qbar)/T' \cong E(\Qbar(t))/E(\Qbar(t))_{\tors}
\]
of lattices and $\Gal(\Qbar/\QQ)$-modules. Now, $\KXperp$ is an $E_8$-lattice with maximal Galois action $W(E_8)$, because $X$ has maximal Galois action on its geometric Picard group by Theorem~\ref{T:main theorem}.  To complete the proof of the theorem, it remains to check that $E(\Qbar(t))_{\tors} = 0$; this is true by~\cite[Theorem 10.4]{Shioda2}
\qed


\end{document}